\documentstyle[12pt,epsf,amssymb,graphicx]{article}

\font\tenmsb=msbm10
\font\sevenmsb=msbm7
\font\fivemsb=msbm5
 
\newfam\msbfam 
\textfont\msbfam=\tenmsb
\scriptfont\msbfam=\sevenmsb
\scriptscriptfont\msbfam=\fivemsb
 
\def\Bbb#1{\fam\msbfam\relax#1}

\newtheorem{thm}{Theorem}[section]
\newtheorem{prop}[thm]{Proposition}
\newtheorem{cor}[thm]{Corollary}
\newtheorem{lem}[thm]{Lemma}
\newtheorem{conj}[thm]{Conjecture}
\newtheorem{exa}[thm]{Example}
\newtheorem{defn}[thm]{Definition}

\newtheorem{rem}[thm]{Remark}
\newtheorem{note}[thm]{Notation}
\newtheorem{alg}[thm]{Algorithm}

\newcommand{\ben}{\begin{enumerate}}
\newcommand{\een}{\end{enumerate}}
\newcommand{\ble}{\begin{lem}}
\newcommand{\ele}{\end{lem}}
\newcommand{\bth}{\begin{thm}}
\renewcommand{\eth}{\end{thm}}
\newcommand{\bpr}{\begin{prop}}
\newcommand{\epr}{\end{prop}}
\newcommand{\bco}{\begin{cor}}
\newcommand{\eco}{\end{cor}}
\newcommand{\bcon}{\begin{conj}}
\newcommand{\econ}{\end{conj}}
\newcommand{\bde}{\begin{defn}}
\newcommand{\ede}{\end{defn}}
\newcommand{\bex}{\begin{exa}}
\newcommand{\eex}{\end{exa}}
\newcommand{\brem}{\begin{rem}}
\newcommand{\erem}{\end{rem}}
\newcommand{\bnot}{\begin{note}}
\newcommand{\enot}{\end{note}}
\newcommand{\balg}{\begin{alg}}
\newcommand{\ealg}{\end{alg}}

\newcommand{\bib}{thebibliography}
\newcommand{\qed}{\square}

\newcommand{\C}{{\Bbb C}}

\newcommand{\PP}{{\Bbb P}}

\begin{document}

\title{Identifying Half-Twists Using Randomized Algorithm Methods} 

\author{S. Kaplan and M. Teicher}


\maketitle

\stepcounter{footnote}\footnotetext{Partially supported by the Emmy Noether Research Institute for
Mathematics, (center of the Minerva foundation of Germany), the Excellency Center ``Group Theoretic 
Methods in the Study of Algebraic Varieties'' of the Israel Science Foundation, and EAGER (EU 
network, HPRN-CT-2009-00099).}
\stepcounter{footnote}\footnotetext{This paper is part of the first author PhD thesis}

\begin{center}
Corresponding address: Shmuel Kaplan

Department of Mathematics and Computer Sciences

Bar Ilan University 

Ramat-Gan, Israel 52900

kaplansh@macs.biu.ac.il

\end{center}

\small{key words: algorithm, random, half-twist, conjugacy, braid group}

\begin{abstract} 
Since the braid group was discovered by E. Artin \cite{Artin}, the question of its
conjugacy problem has been solved by Garside \cite{GAR} and Birman, Ko and Lee \cite{NEW}. However, 
the solutions given thus far are difficult to compute with a computer, since the number of 
operations needed is extremely large. Meanwhile, random algorithms used to solve difficult problems 
such as primality of a number were developed, and the random practical methods have become an
important tool. We give a random algorithm, along with a conjecture of how to improve its
convergence speed, in order to identify elements in the braid group, which are 
conjugated to its generators (say $\sigma _1^k$) for a given power $k$. These elements of the braid 
group, the half-twists, are important in themselves, as they are the key players in some 
geometrical and algebraical methods, the building blocks of quasipositive braids and they construct 
endless sets of generators for the group.

\end{abstract}

\section{Introduction} 

Let $B_n$ be the braid group on $n$ strings. 
The conjugacy problem in $B_n$ is difficult and was addressed in several cases
in the past. The computation of a solution is still not accessible. Although
a solution was proposed (\cite{GAR},\cite{POSBR}), the running time of a computerized
program based on these algorithms is extremely long.

In the sequel we will describe an algorithm that solves a partial problem of the conjugacy
problem. Our algorithm will make it possible to identify whether for a given braid $w$ there exists 
an integer $k$ such that $w$ is conjugated to $\sigma _1^k$. Therefore we actually identify some 
special conjugacy classes of the braid group.

Our algorithm is based on a random technique, and has the property that in any case that
it returns $true$ (meaning that the input element $w$ of the braid group is conjugated to a
generator of the group in some power), it will also return to which generator $\sigma _i$ the 
element is conjugated, and what is $q$ such that $q^{-1}wq=\sigma _i^k$. Although we do not have 
any estimations for the probabilitiy of an erroneous return value, we did a large number of 
experiments that gave us information as to how well the algorithm converges.

Elements in the conjugacy class of one of Artin generators are called half-twists. If one uses 
braid monodromy in order to classify geometrical hypersurfaces up do deformation 
\cite{BGTI}\cite{KUL}, this will result with half-twists. Therefore, by using this algorithm it is 
possilbe to verify braid monodromy computations. Moreover, identifying the elements of this 
conjugacy class has implications to the research of quasipositive braids.

We start in Section 2 by giving some braid group preliminaries. In Section 3 we give a complete 
description of the
random method for identifying half-twists along with full proofs of its correctness and complexity. 
Section 4 is dedicated to the presentation of the experiments and benchmarks that
were done in order to understand the capabilities of the random algorithm. Finally, Section 5 is 
dedicated to closing arguments.

\section{Braid group preliminaries}

\subsection{E. Artin's definition of the braid group}

\bde
\underline{Artin's braid group $B_n$} is the group generated by $\{\sigma _1,...,\sigma _{n-1}\}$
subject to the relations 
\ben
\item $\sigma _i\sigma _j=\sigma _j\sigma _i$ where $|i-j| \geq 2$ 
\item $\sigma _i\sigma _{i+1}\sigma _i=\sigma _{i+1}\sigma _i\sigma _{i+1}$  for all $i=1,...,n-2$
\een
\ede

This algebraic definition can be looked at from a geometrical point of view, by associating to 
every generator
of the braid group $\sigma _i$ a tie between $n$ strings going monotonically from top to bottom, 
such
that we switch by a positive rotation between the two adjacent pair of strings $i$ and $i+1$. 
This means that $\sigma _i$ corresponds to the geometrical element described in the following 
figure:

\begin{figure}[h]
\centering
\includegraphics[width=2cm]{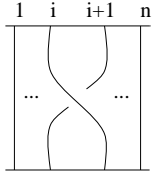}
\caption{The geometrical braid associated with $\sigma _i$}
\end{figure}

The operation for the geometrical group is the concatenation of two geometrical 
sets of strings. 

\bex
The geometrical braid that corresponds to $\sigma _1 \sigma _2^{-1} \sigma _1 \sigma _3$ 
is presented in the following figure:

\begin{figure}[h]
\centering
\includegraphics[width=2cm]{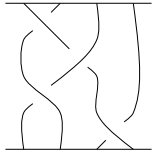}
\caption{The geometrical braid $\sigma _1 \sigma _2^{-1} \sigma _1 \sigma _3$}
\end{figure}
\eex

\subsection{Half-twists in the braid group}

We are now going to describe what half-twists are, and what is their geometrical
interpretation.

\bde
Let $H$ be the conjugacy class of $\sigma _1$, (i.e. $H=\{q^{-1} \sigma_1 q : q \in
B_n$\}). We call $H$ \underline{the set of half-twists in $B_n$}, and we call an element $\beta 
\in H$ a \underline{half-twist}.
\ede

Recall that half-twists have a geometrical interpretation. One can look at the braid group as the 
mapping class group of an $n$-punctured disk. The half-twists then correspond to geometric 
half-twists around an embedded arc that connects two punctures, and do not intersect itself or any 
other puncture. Using this way to look at the half-twists it is easy to see that they occupy a full 
conjugacy class of the braid group, and that they all conjugate to Artin's generators of the group.

Our main goal now is to describe the algorithm, but before we can do that we need some definitions.

\bde
Let $w \in B_n$ be a braid. Then it is clear that $w=\sigma _{i_1}^{e_1} \cdot ...\cdot \sigma
_{i_l}^{e_l}$ for some sequence of generators, where $i_1,...,i_l \in \{1,...,n-1\}$ and 
$e_1,...,e_l \in \{1,-1\}$. We will call such a presentation of $w$ a \underline{braid
word}, and $\sigma _{i_k}^{e_k}$ will be called the \underline{$k^{th}$ letter of the word $w$}. 
We will call $l$ the \underline{length} of the braid word.
\ede

We will distinguish between two relations on braid words.

\bde
Let $w_1$ and $w_2$ be two braid words. We will say that $w_1=w_2$ if they represent the
same element of the braid group.
\ede

\bde
Let $w_1$ and $w_2$ be two braid words. We will say that $w_1 \equiv w_2$ if $w_1$ and
$w_2$ are identical letter by letter.\
\ede

\bde
A \underline{positive braid} is an element of $B_n$ which can be written as a word in positive 
powers of the
generators $\{\sigma _i\}$, without the use of the inverse elements $\sigma _i^{-1}$. We denote
this subsemigroup $B_n^+$
\ede

Next, we are going to recall some of the algorithms that were developed by A. Jacquemard, in
\cite{EFF}. These algorithms perform manipulation over braid words, and we
will use a combination of them in our random method for determining if an element of the braid
group is a half-twist or not. 

The $LetterExtractLeft$ algorithm enables us to determine if it is possible to write a certain
positive braid word in such a way that a given generator is its first letter.

\ble \cite{EFF} 
If $l$ is the length of the word $w$, then the $LetterExtractLeft$ procedure complexity is 
$O(l^2)$.
\ele

\brem
The $LetterExtractLeft$ procedure can be altered easily to become the $LetterExtractRight$ 
procedure.
\erem

The $WordExtractLeft$ algorithm tries to drive a given positive braid word to the left of a 
positive braid word. It determines if it is possible to write the positive braid word $w$ in such a 
way that the
left part of $w$ will be a given positive braid word.

We make use of this algorithm as a function in the following way: \\
$WordExtractLeft(w,w')$, its inputs are a positive braid word $w$ from which we are going
to extract the positive braid word $w'$ to the left. If it is not possible to write
$w=w'w''$ for some positive braid word $w''$, the function will return {\it false}. If it is 
possible to write $w=w'w''$ for some positive braid word $w''$, the return value will be $w'w''$ 
(i.e. the braid $w$ written as a positive word where $w'$ is its prefix).

\ble \cite{EFF} 
If $l_1$ is the length of $w$ and $l_2$ is the length of $w'$, then the complexity of the 
algorithm is $O(l_2 \cdot l_1^2)$.
\ele

\brem
The $WordExtractLeft$ procedure can be altered easily to become the $WordExtractRight$ procedure.
\erem

The $Normalize$ algorithm transforms a given braid word $w$ into a normal form
$w=\Delta _n^{-r}w'$, where $\Delta _n$ is the fundamental braid word, $r$ is minimal 
and $w'$ is a positive braid word given in its lowest lexicographical order possible (This is 
Jacquemard's solution to the braid word problem and it is cubic in the length of the braid word).

\section{The randomized algorithm}
In this section we will describe our random algorithm for identifying half-twists in any
power. First we will describe what a random algorithm is and how it works 

\subsection{Randomized algorithm}
Let $D$ be the set of all possible inputs to an algorithm $\cal{A}$, which is supposed to compute a 
function $A$. Consider the situation where $D$ is divided into equivalence classes by the relation 
$\sim$. Suppose that $A:D/\sim \to \{0,1\}$ is well defined mathematically.
If $\cal{A}$ is an algorithm that will return the wrong answer on some inputs given from $D$,
but on others will return the right answer, and if we have a method to make sure that the answer 
is correct or not in some cases, this makes a solid foundation for using it as a
random algorithm.

In our case, $D$ is the set of braid words $w$ represented by the generators in both
positive and negative powers. The equivalence relation $\sim$ is given by $w_1 \sim w_2$ if
they are conjugated. This actually means that the half-twists occupy a full equivalence class of 
$\sim$, and that two half-twists in the same power share the same equivalence class of $\sim$. 
 
The function $A$ returns {\it true} if $w$ is conjugated to a half-twist in some positive power and
{\it false} otherwise. We will present an algorithm for solving the problem $\cal{A}$, but
unfortunately we will only be able to be convinced of its result if its answer is {\it
true}, or in some cases of a {\it false} answer, which we will call a \underline{{\it genuine} {\it 
false}}. In some cases, although
its input $w$ is a half-twist in some positive power, our algorithm will return {\it false}. 

As we will state in the next sections the probability of such an error is low;
therefore it will be possible to iterate the algorithm on different elements of the
conjugacy class of the input, resulting in a substantial reduction in the probability of the
error. It is possible to reduce the probability of an error as much as possible, simply by 
increasing the number of iterations.

One very important and known example for a random algorithm is the algorithm for checking whether a
natural number is prime or not. A randomized method developed by Miller \cite{MILLER} and
Rabin \cite{RABIN}, using iteration of 
checking 'witnesses' to the primality of the given number, results in a probability of
an error which can be reduced as much as one wants. Calculation of large prime
numbers, especially with a connection to encryption, is done today using this method.

\subsection{Check for conjugacy to $\sigma _i$}
In this section we will describe the algorithm for checking whether or not a given braid word $w$
is a half-twist in some power. The algorithm is based on two functions that will be
described in detail in the subsections below.

The idea of the algorithm is based on the fact that if $w$ is a half-twist, then there exists a 
braid
word $q$ such that $q^{-1}wq=\sigma _i$, for any $\sigma _i$ generator of the braid group.
First we find the power of the alleged half-twist by summing up the powers of the
generators in $w$; this power will be denoted by $k=deg(w)$. Then we will try to manipulate the
braid word $w$ to be written as $w=q^{-1} \sigma _i^k q$; this will be done
using the algorithms for braid word manipulation given above.

If the result of the word manipulation is {\it true} (meaning that we achieved the form $w=q^{-1} 
\sigma
_i^k q$,), then we return {\it true}, since the braid word is obviously
a half-twist. Moreover, we can return $q$ and $\sigma _i$ which are the way the word $w$
is conjugated to a generator of the braid group in the $k^{th}$ power. 

On the contrary, if the result of the algorithm is {\it false}, still it could be possible that the 
input braid word $w$ is conjugated to a generator in the $k^{th}$ power. 

However, this problematic situation mentioned above is quite rare, allowing us the possibility
to use the random iteration method. In this case we choose a random braid word in the length of $w$ 
and conjugate $w$ by it, resulting in a new braid word $w'$ which
is conjugated to $\sigma _i^k$ if and only if $w$ is conjugated to $\sigma _i^k$. Then we
try the algorithm again on $w'$.

We conjecture, by looking at the data from the experiments that we did, that the probability of an 
error in the result does not change significantly when we build the new element $w'$ in the 
conjugacy class of the given word $w$. Therefore, although we do not know how to prove that 
rigorously, we believe that if the probability of an error in the result is smaller than $p<1$ 
then, by iterating $n$ times, one can get to a certainty of nearly $1-p^n$ that the given braid
word $w$ is not a half-twist. 

In the performance and benchmarks section we will give our estimations to the
error rate.

\subsection{The Is half twist procedure}

First we give the algorithm that tries to manipulate the input braid word $w$ into
$w=q^{-1} \sigma _i^kq$, where $k$ is the sum of powers of the braid word $w$.

But, before we can do that we need to introduce a notation.

\bnot
Let $w$ be a braid word. We denote by $w_{i,j}$ the part of $w$ that starts at the
$i^{th}$ letter and ends at the $j^{th}$ letter of $w$. If $j<i$ then $w_{i,j}$ denotes the empty 
word.
\enot

\balg{IsHalfTwist}

\noindent
${\bf Input:}$ $w$ - non empty braid word in its normal form. $k$ - the power of the half-twist.

\noindent
${\bf Output:}$ true - if the element is a half-twist to the $k^{th}$ power, false if the element 
is not a half-twist to the $k^{th}$ power. $c$ - the generator for which we found the
conjugacy, $q$ - a braid word which conjugate $w$ to $c$ (i.e., $q^{-1}wq=c^k$).

\parskip0pt

\noindent 
{\bf IsHalfTwist}($w,k$)

$iPos \leftarrow$ the position of the first positive letter in $w$

$l \leftarrow$ the length of $w$

{\bf for} $i \leftarrow 1$ {\bf to} (the number of strings in $B_n$) $-$ $1=n-1$ {\bf do}

$\qquad$ {\bf for} $t \leftarrow 0$ {\bf to} $k$ {\bf do}

$\qquad$ $\qquad$ $pLeft \leftarrow \sigma _i^t$

$\qquad$ $\qquad$ $pRight \leftarrow \sigma _i^{k-t}$

$\qquad$ $\qquad$ {\bf for} $p \leftarrow iPos-1$ {\bf to} $l$ {\bf do}

$\qquad$ $\qquad$ $\qquad$ {\bf if} $WordExtractRight(w_{iPos,p},pLeft)=false$ {\bf then}

$\qquad$ $\qquad$ $\qquad$ $\qquad$ {\bf continue}

$\qquad$ $\qquad$ $\qquad$ {\bf if} $WordExtractLeft(w_{p+1,l},pRight)=false$ {\bf then}

$\qquad$ $\qquad$ $\qquad$ $\qquad$ {\bf continue}

$\qquad$ $\qquad$ $\qquad$ $Test \leftarrow$ $w_{0,p-t}w_{p+k-t+1,l}$

$\qquad$ $\qquad$ $\qquad$ $q \leftarrow$ $(w_{p+k-t+1,l})^{-1}$

$\qquad$ $\qquad$ $\qquad$ {\bf if} $Test=Id$ {\bf then}

$\qquad$ $\qquad$ $\qquad$ $\qquad$ {\bf return} $true,q,\sigma _i$

{\bf return} false

\ealg
\parskip3pt

\bpr
Given a braid word $w$ if the algorithm returns {\it true}, then the braid word $w$ is a
half-twist in the $k^{th}$ power.
\epr

\noindent
{\it Proof:} 
The algorithm tries to write $w$ as $q\sigma _i^kq'$ for each
position in the positive part of $w$ and each generator possible for $\sigma _i^k$ $(i=1,...,n-1)$. 
Then, it checks if $q=q'^{-1}$. If this happens, the algorithm will return {\it true}. But this is 
exactly the case when $w$ can be written as a
conjugacy to a half-twist in the $k^{th}$ power, meaning that $w$ is a half-twist in the
$k^{th}$ power. Moreover, if $w$ is not conjugated to $\sigma _i^k$, then there is no way
to write $w$ as $q\sigma _i^kq'$ where $q=q'^{-1}$, and since the only possible return of {\it 
true} is after such
writing occur, the algorithm will not return {\it true} if $w$ is not a half-twist to the
$k^{th}$ power.
\hfill $\qed$

We give two examples to illustrate why the algorithm can return {\it false} although its result 
should be {\it true}.

\bex
Look at the braid word $\sigma _1 ^{-1} \sigma _2 \sigma _1 \sigma  _2 \sigma _1^{-1}$. There is no 
partition of the word as the algorithm does that will change it into $q^{-1}\sigma _i q$ for some 
braid word $q$, although $\sigma _1 ^{-1} \sigma _2 \sigma _1 \sigma  _2 \sigma _1^{-1}=\sigma _2$. 
\eex

Actually, since the algorithm works only on words in Garside normal form, this braid word does not 
cause any incorrect {\it false}. The next example gives us a braid word which does cause the 
algorithm to return an incorrect {\it false}.

\bex
The following braid word in $B_5$ is a half-twist (conjugated to $\sigma _1$) but the algorithm
will return {\it false} while processing on it:

$\sigma _4^{-1}\sigma _1\sigma _3^{-1}\sigma
_1^{-1}\sigma _3^{-1}\sigma _4\sigma _4\sigma _1^{-1}\sigma _1\sigma _2^{-1}\sigma
_1\sigma _2^{-1}\sigma _1\sigma _2\sigma _1^{-1}\sigma _2\sigma _1^{-1}\sigma _1\sigma
_4^{-1}\sigma _4^{-1}\sigma _3 \sigma _1\sigma _3\sigma _1^{-1}\sigma _4$.
\eex

We now give an explanation why the algorithm might return incorrect {\it false} results. Let 
$w=q^{-1}\sigma _i^kq$ be a half-twist. Therefore, there is a finite sequence of braid relation 
actions on $w$ that will transform it into $w' \equiv q^{-1} \sigma _i^k q$. We keep track of the 
position in the word of each of the letters $\sigma _i$ that belongs to the original $\sigma _i^k$ 
of $w'$. Suppose that we only use the first braid relation ($\sigma _i\sigma _j=\sigma _j\sigma _i$ 
where $|i-j|>1$). It is clear that if we eliminate the letters from the original $\sigma _i^{k}$ in 
any stage, this results in a word that equals to $Id \in B_n$. Therefore, activation of the first 
braid relation will never cause a braid word to result in an incorrect result. This leaves only the 
second rule as a possible cause, and so we have the following lemma:

\ble
Let $w=q^{-1}\sigma _1^k q$ be a braid word such that the $IsHalfTwist(w,k)$ function returns {\it 
false}. Then, in the process of transforming the braid word from $q^{-1}wq$ into $w$ we must use 
the second braid relation $\sigma _i \sigma _{i+1} \sigma _i=\sigma _{i+1} \sigma _i \sigma 
_{i+1}$.
\ele
\hfill $\qed$

\bpr
The complexity of the $IsHalfTwist$ algorithm is bounded by $O(n^2 \cdot log(n) \cdot k^2 \cdot 
l^3)$, where $n$ is
the number of strings, $k=deg(w)$ and $l=len(w)$.
\epr

\noindent
{\it Proof:}
The proof involves the following facts: $(a)$ We go over all possible
generators $\sigma _i$ $1 \leq i \leq n-1$. $(b)$ We go over all possibilities to divide the
$\sigma _i^k$ into two subwords. $(c)$ We go over all positions $p$ in the positive braid
part of the braid word $w$. The total gives us $O(n \cdot k \cdot l)$ times which we call the 
$WordExtractLeft$ and $WordExtractRight$ procedures. The Complexity of the
$WordExtractRight$ procedure in our case is given by $O(k' \cdot l'^2)$,
where $k'$ is the power of the left part of $\sigma _i^k$ and $l'=p-iPos$ at each loop. The
complexity of the $WordExtractLeft$ procedure in our case is given by $O(k'' \cdot l''^2)$
where $k''$ is the power of the right part of $\sigma _i^k$ and $l''=l-p$. Together we
have that $k'+k''=k$ and $l'+l''=l-iPos$, and therefore, the number of operations we make in the 
two procedures is bounded by $O(k \cdot l^2)$. Combining that result with the above, we get a total 
of $O(n \cdot k^2 \cdot l^3)$. To finish the proof, we need to observe that the check whether 
$Test=Id$ takes $O(n \cdot log(n) \cdot l^2)$ (see for example \cite{NEW} \cite{Deh1} \cite{Deh2} 
for a fast solution for the braid word problem). This yields the total of $O(n^2 \cdot log(n) \cdot 
k^2 \cdot l^3)$ as the total complexity bound for the $IsHalfTwist$ algorithm.

\hfill $\qed$

\subsection{Some ways of improving the running time}

There are some ways to make the algorithm run faster. The first two shortcuts involve
keeping track of what exactly happens with the strings of the braid. So, we will enumerate
the strings by $1,...,n$.

\bde
We call the number of times that two strings $i \neq j$ cross one another (counted with
the positivity induced by the sign of the generator denoting the switch) their
\underline{crossing number}, and we denote it by $cr(i,j)$.
\ede

\ble
If $w$ is a half-twist to the $k^{th}$ power, and $k$ is even, there must be only one pair of 
strings with crossing number $cr(i,j) \neq 0$, and then $cr(i,j)=k$.
\ele

\noindent
{\it Proof:} It follows from the fact that using braid relations do not change the
crossing numbers of pairs of strings in $w$, and from the geometrical interpretation of
the half-twist elements.
\hfill $\qed$

Therefore, if $k$ is even we can begin our algorithm by counting the crossing numbers of each pair 
of strings. If $cr(i,j) \neq 0$ for more than one pair, or $cr(i,j)=0$ for all $i \neq j$,
we can return $false$ which is genuine (i.e. the probability of this $false$ result to be
correct is $1$).

Moreover, if $w$ is a half-twist and we eliminate from $w$ the two strings whose crossing number is 
not $0$, resulting with $w' \in B_{n-2}$, then $w'=Id \in B_{n-2}$. Therefore, if $w' \neq Id \in 
B_{n-2}$, we can return a {\it genuine} {\it false}.

\bde
Assign to each letter in $w$ a pair of numbers which represents the numbers of the
strings that switch because of this generator; we call this pair the
\underline{switching numbers} of the letter.
\ede

By keeping track of all the switching numbers of the letters of $w$ from the beginning to
the end (no need to defer between positive or negative power of the letter in $w$), we result in an 
easy way to compute the permutation of the strings at each part of the word
$w$. After looking at all the letters, we get the permutation of the strings resulted by
$w$. 

\ble
Let $w$ be a braid word representing a half-twist to the $k^{th}$ power. Then the
permutation induced by $w$ should be $Id$ if $k$ is even or a transposition of exactly two
strings if $k$ is odd.
\ele

\noindent
{\it Proof:} Immediately follows from the fact that a half-twist always changes the
position of two strings along a path between the strings (see the geometrical interpretation of the 
half-twists). 
\hfill $\qed$

Now, we can check if the permutation is consistent with $k$ which
is the sum of the powers of the letters in $w$. If there is inconsistency, we can return a {\it 
genuine} $false$. 

Moreover, if $k$ is odd, $w$ is a half-twist and we eliminate from $w$ the two strings that 
permute, we result with $w' = Id \in B_{n-2}$. Therefore, if $w' \neq Id \in B_{n-2}$ we can return 
a {\it genuine false}.

After doing these two steps, we know what two strings are switching position.
Therefore, checking conjugacy to another pair of switching strings is irrelevant, and can
be skipped. This means that after we have written our braid word as $w=q\sigma _i^kq'$, we
can check if $\sigma _i^k$ induces the switch of the two strings we found at the
beginning. If not, there is no need to check if $q=q'^{-1}$. 

It seems that we need to compute what are the switching pair of strings at the position
$\sigma _i^k$ each time we manage to write our braid as $w=q\sigma _i^kq'$ but this is
not true. One can easily create the switching pairs for all the letters by one pass over the given 
braid word $w$, and then maintain easily the switching pair data by following the next
rules:

\ble
In the braid manipulation process, one of the following must occur:

\ben
\item
If we use the first braid rule $\sigma _i \sigma _j=\sigma _j \sigma _i$, then we must
switch the two appropriate switching pairs.
\item
If we use the second braid rule $\sigma _i \sigma _{i+1} \sigma _i=\sigma _{i+1} \sigma _i
\sigma _{i+1}$, then we must switch the two switching pairs assigned to the two ends of
the sequence.
\een
\ele

\noindent
{\it Proof:}
To begin with, it is clear that no other relation is used in the braid manipulation process;
therefore we need only to proof that these are the changes in the switching pairs that are
induced by the two relations.

For the first rule, since $\sigma _i \sigma _j=\sigma _j \sigma _i$ we know that
$|i-j|>1$. This means that the two switching pairs are distinct. Therefore, the switching
pair associated with $\sigma _i$ and $\sigma _j$ continue with their letters, and they do
not collide. 

For the second rule, we know that the switching pairs must be of the form
$(a,b),(a,c),(b,c)$. Therefore, the change in the letters induces the change of the switching pairs 
to $(b,c),(a,c),(a,b)$ (see the diagram below).

\begin{figure}[h]
\centering
\includegraphics[width=6cm]{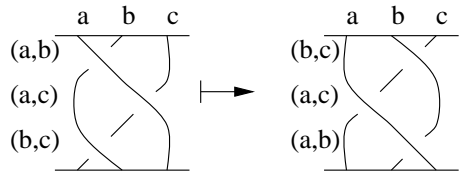}
\caption{The second rule of switching pairs}
\end{figure}

\hfill $\qed$

The last improvement in the running time of the algorithm we are going to present is a
consequence of a $false$ return value of the $LetterExtractLeft$ or the $LetterExtractRight$
functions. Suppose that we are trying to 
extract to the right the letter $\sigma _i$ at position $p$, and we do not succeed.
This means that there is no chance of extracting the letter $\sigma _i$ to the right until we
encounter the next $\sigma _i$ in the braid word. Therefore, any attempt to extract to the
right the letter $\sigma _i$ for positions from $p$ to the position of the next $\sigma
_i$ letter in $w$ will fail. Therefore, we can skip these tries without
affecting the result of the algorithm. Note, that the same is true for the
$LetterExtractLeft$, but this time we should check what is the leftmost position that we
could extract the letter $\sigma _i$ to the left and move $p$ immediately to this position. Note 
that due to the way that the $LetterExtractLeft$ function works, this
leftmost position is given instantly after the first unsuccessful 
$LetterExtractLeft$ call. 

This enables us to skip a large amount of tries to extract the letters into positions that
we know already we can't extract to. By keeping track of what is the leftmost or rightmost
position we can extract each letter, we can make the algorithm run faster.

We want to clarify that although in some cases it is possible to write the word with the letter 
that was extracted to the left (or to the right) at a position which is left to its position as 
returned by the $LetterExtractLeft$ algorithm, since the $IsHalfTwist$ procedure does not consider 
these situations, skipping extraction to these positions does not change the results.

\subsection{The Test Random Half-Twist algorithm}

In this section, we will give the random core of the algorithm, the $TestRandomHalfTwist$
procedure. This procedure gets the braid word and tries successively to check if it is a
half-twist. The procedure will terminate in one
of three conditions. The first is when $IsHalfTwist$ returns {\it true}, this happens when no
question about the correctness of the answer exists, therefore one can return the computed result. 
The second is when we encounter a {\it genuine} {\it false}; here again there is no doubt of
the correctness of the result. The third is when we have exceeded some predefined maximum value of 
tries. 
In this situation we say that the given braid word $w$ is not conjugated to a half-twist,
with the restriction that there could be an error which its probability is less than $p$ (where 
$p$ will be discussed later).

\balg{Test Random Half-Twist procedure}

\noindent
${\bf Input:}$ $w$ - non empty braid word in its normal form. 

\noindent
${\bf Output:}$ true - if the element is a half-twist to the $k^{th}$ power. false - if the element 
is not a half-twist to the $k^{th}$ power, or a genuine false if possible. $c$ - the generator for 
which we found the
conjugacy. $q$ - a braid word which conjugate $w$ to $c$ (i.e., $q^{-1}wq=c^k$).

\parskip0pt

\noindent 
{\bf TestRandomHalfTwist}($w$)

$r \leftarrow$ $Id$

$r1 \leftarrow w$

$i \leftarrow 0$

$k \leftarrow deg(w)$

{\bf do}

$\qquad$ $i \leftarrow i+1$

$\qquad$ {\bf if} $IsHalfTwist(r1,k)=true$ {\bf then}

$\qquad$ $\qquad$ $q \leftarrow r \cdot q$

$\qquad$ $\qquad$ $c \leftarrow$ the letter where $IsHalfTwist$ found the conjugacy

$\qquad$ $\qquad$ {\bf return} true

$\qquad$ {\bf else}

$\qquad$ $\qquad$ {\bf if} the return value was a genuine false {\bf then}

$\qquad$ $\qquad$ $\qquad$ {\bf return} a genuine false.

$\qquad$ $\qquad$ {\bf else}

$\qquad$ $\qquad$ $\qquad$ $r \leftarrow$ a random braid word of size $length(w)$

$\qquad$ $\qquad$ $\qquad$ $r1 \leftarrow Normalize(r^{-1} \cdot w \cdot r)$

{\bf while} $i<$ maximum number of tries

{\bf return} false (not genuine)

\ealg
\parskip3pt

\bpr
Given a braid word $w$ when the algorithm returns {\it true}, the result $q$ is the
conjugacy word for which $q^{-1}wq=c^k$. Moreover, on a half-twist the probability of an incorrect 
false return value is $p' \cdot p^{i-1}$, where $p$ is the probability of an error for the result 
of $IsHalfTwist$ on conjugated words to $w$, $p'$ is the probability of an error for the result of 
$IsHalfTwist$ on $w$ itself and $i$ is the number of tries.
\epr

\noindent
{\it Proof:}
The first part is obvious, since the only reason that this procedure returns {\it true} is when
$IsHalfTwist$ returns {\it true}. 

Now, we have to show that the returned conjugacy word $q$
satisfies that $q^{-1}wq=c^k$. We have two cases. The
first case is when the algorithm returns {\it true} in the first pass of the 'do' loop. In this
case, the braid word that was entered to the $IsHalfTwist$ procedure is $w$, and the returned $q$
is exactly as the $q$ given by $IsHalfTwist$; therefore, this is the correct $q$. In the
second case, we have iterated more than once over the 'do' loop. At each iteration we
have conjugated by $r$ the braid word $w$ resulting with $r1=r^{-1}wr$. Therefore, the $q$
returned by the $IsHalfTwist$ procedure satisfies that $q^{-1} \cdot r1 \cdot q=c^k$. But
$r1=r^{-1}wr$; therefore, $q^{-1}r^{-1}wrq=c^k$, as the algorithm returns.

For the second part of the proposition suppose that the probability to receive an error by
activating $IsHalfTwist$ on $w$ is $p'$, and that the probability to receive an error by activating 
$IsHalfTwist$ on the random conjugates of $w$ is $p$. Consider the fact that if the algorithm has 
returned {\it false} this means that the 'do' was executed $i$ times. Therefore, the probability of 
an erroneous result is $p'p^{i-1}$.
\hfill $\qed$

\bpr
The complexity of the $TestRandomHalfTwist$ procedure is $O(n^2 \cdot log(n) \cdot k^2 \cdot l^3)$.
\epr

\noindent
{\it Proof:}
Since the number of times we activate the procedure $IsHalfTwist$ is constant, and since all
the other procedures we do in the 'do' loop take constant time or less then the time for
the $IsHalfTwist$ procedure, this results in a constant times the complexity of the
$IsHalfTwist$ procedure which is $O(n^2 \cdot log (n) \cdot k^2 \cdot l^3)$.

\hfill $\qed$

\section{Performance}

In this section we will give the probability estimations for the success rate for the
$IsHalfTwist$ procedure. We will also state that one may check $w^2$ instead of
$w$, reducing dramatically the probability for an incorrect {\it false}.

\subsection{Probability estimations}

Since the combinatoric computations of the braid group that are needed in order to give a result on 
the
probability of an incorrect {\it false} return value are yet unknown, we have made
numerous experiments using a computer program. The results are encouraging since it
looks as if the probability to get an incorrect {\it false} return value from the $IsHalfTwist$ 
procedure is low.

In our experiments we shuffled random half-twists and transformed them into normal form.
Then, we activated the $IsHalfTwist$ once, counting for each length of word the number
of times we tried the algorithm and the number of times that the algorithm resulted with
{\it true}. We did this in the groups $B_5$ and $B_8$ for powers of the
half-twists between $1$ and $5$. The results are summarized in the tables below. Each line 
represents the
tests on words of length in an interval of 100 letters. The two numbers in each Power column 
represent the number of words shuffled and the probability of success of the algorithm.

\begin{center} 

\begin{table}
\caption{Success probability of the algorithm on $B_5$}
\centering
\begin{tabular}{|c|c|c|c|c|c|}
\hline
\hline  {length}    & {Power 1}     & {Power 2}      & {Power 3}  & {Power 4}  & {Power 5}  \\
\hline  {100}       & {535 , 0.985} & {673 , 1}      & {1149 , 1} & {1139 , 1} & {1196 , 1} \\ 
\hline  {200}       & {550 , 0.829} & {684 , 1}      & {1126 , 1} & {1208 , 1} & {1131 , 1} \\
\hline  {300}       & {608 , 0.781} & {671 , 0.997}  & {1087 , 1} & {1046 , 1} & {1118 , 1} \\
\hline  {400}       & {553 , 0.772} & {555 , 0.994}  & {941 , 1}  & {905 , 1}  & {881 , 1}  \\ 
\hline  {500}       & {492 , 0.770} & {311 , 0.993}  & {558 , 1}  & {542 , 1}  & {532 , 1}  \\
\hline  {600}       & {452 , 0.803} & {92 , 1}       & {127 , 1}  & {139 , 1}  & {116 , 1}  \\
\hline  {700}       & {297 , 0.771} & {14 , 1}       & {12 , 1}   & {21 , 1}   & {26 , 1}   \\
\hline  {800}       & {109 , 0.779} & {-}            & {-}        & {-}        & {-}        \\
\hline  {900}       & {44 , 0.701}  & {-}            & {-}        & {-}        & {-}        \\
\hline  {1000}      & {-}           & {-}            & {-}        & {-}        & {-}        \\
\hline
\end{tabular}
\end{table}

\begin{table}
\caption{Success probability of the algorithm on $B_8$}
\begin{tabular}{|c|c|c|c|c|c|} 
\hline  {length} & {Power 1}     & {Power 2}     & {Power 3} & {Power 4} & {Power 5} \\
\hline  {100}    & {147 , 1}     & {135 , 1}     & {141 , 1} & {762 , 1} & {488 , 1} \\ 
\hline  {200}    & {125 , 0.944} & {142 , 1}     & {148 , 1} & {685 , 1} & {476 , 1} \\
\hline  {300}    & {135 , 0.755} & {124 , 1}     & {138 , 1} & {634 , 1} & {427 , 1} \\
\hline  {400}    & {136 , 0.669} & {135 , 1}     & {130 , 1} & {678 , 1} & {443 , 1} \\ 
\hline  {500}    & {64 . 0.625}  & {58 , 1}      & {59 , 1}  & {297 , 1} & {195 , 1} \\
\hline  {600}    & {64 , 0.577}  & {118 , 0.991} & {100 , 1} & {556 , 1} & {410 , 1} \\
\hline  {700}    & {107 , 0.598} & {100 , 1}     & {107 , 1} & {511 , 1} & {417 , 1} \\
\hline  {800}    & {91 , 0.466}  & {89 , 0.977}  & {77 , 1}  & {411 , 1} & {425 , 1} \\
\hline  {900}    & {55 , 0.527}  & {58 , 0.931}  & {54 , 1}  & {142 , 1} & {397 , 1} \\
\hline  {1000}   & {15  0.733}   & {22 , 1}      & {17 , 1}  & {189 , 1} & {188 , 1} \\
\hline  {1100}   & {17 , 0.470}  & {11 , 1}      & {18 , 1}  & {79 , 1}  & {374 , 1} \\
\hline  {1200}   & {-}           & {-}           & {6 , 1}   & {44 , 1}  & {311 , 1} \\
\hline  {1300}   & {-}           & {-}           & {4 , 1}   & {9 , 1}   & {190 , 1} \\
\hline  {1400}   & {-}           & {-}           & {-}       & {1 , 1}   & {73 , 1}  \\
\hline  {1500}   & {-}           & {-}           & {-}       & {2 , 1}   & {103 , 1} \\
\hline  {1600}   & {-}           & {-}           & {1 , 1}   & {-}       & {45 , 1}  \\
\hline  {1700}   & {-}           & {-}           & {-}       & {-}       & {24 , 1}  \\
\hline  {1800}   & {-}           & {-}           & {-}       & {-}       & {10 , 1}  \\
\hline  {1900}   & {-}           & {-}           & {-}       & {-}       & {1 , 1}   \\
\hline  {2000}   & {-}           & {-}           & {-}       & {-}       & {3 , 1}   \\
\hline           
\end{tabular}
\end{table}
\end{center} 

One may suggest that instead of using conjugation by a random braid word, it may be sufficient
to simply increase the number of strings, in order to get a correct result in case of a
mistaken {\it false}. Unfortunately, experiments show that this is not true. If a given braid
word $w$ causes an incorrect {\it false} return value in the $IsHalfTwist$ procedure, then the 
result
is independent of the number of the strings in $B_n$, Therefore, trying to embed $B_n$ in
$B_m$ where $n<m$ will not solve the incorrect {\it false} problem.

\subsection{Square method}

During the check of the different benchmarks, we found that the $IsHalfTwist$ algorithm
almost never returns an incorrect {\it false} result on powers greater than $1$. In addition to the 
results summarized in the tables above, we have tested $13$ half-twists to the $5^{th}$ power in 
$B_8$ with size in normal form of more than 4500 letters, resulting in the return value of {\it 
true} each time. This led us to check the next situation using a computer and we reached the 
following conjecture, which we discuss in \cite{HT}.

\bcon
A braid word $w$ is conjugated to the generator $\sigma _i$ if and only if
$w^2$ is conjugated to $\sigma _i^2$.
\econ

One direction of the proposition is obviously true, since if $w$ is conjugated to the generator
$\sigma _i$ then it is obvious that $w^2$ is conjugated to $\sigma _i^2$.
The surprising part may be the other direction. 

The meaning of this is that it may be possible to reduce the $IsHalfTwist$ problem for general $k$ 
into larger numbers. Given a braid word $w$ with $deg(w)=1$, activate $IsHalfTwist$ on $w^2$ and 
return its result.

Since the probability of an erroneous result in large powers is much lower, we get a better
solution.

Note, that a {\it false} result returned by checking powers of $w$ implies that $w$ is not a 
half-twist in any power.

\section{Conclusions}

We would like to point out that the algorithms given for the $LetterExtractLeft$ and
$WordExtractLeft$, as well as for the $Normalize$, are based on the natural set of generators
for the braid group. It is already known (see for example \cite{BAND}, \cite{NEW}) that if one
changes the set of generators it is possible to improve even more the algorithms complexity. This 
improvement may reflect immediately on the complexity of the random method given
here, if we will find a way to extend the results of \cite{EFF} to work on the new sets of 
generators for $B_n$.

The above algorithm enables us in a fast way to decide whether a given braid word is
conjugated to any generator of the braid group in any power. This makes it possible to
identify these unique and important elements of the braid group, which are the factors in 
a factorization yielded from a braid monodromy. For example, if we use this 
method, it is possible to check whether or not a given factorization of $\Delta ^2$ can be resulted 
from a braid monodromy.

Another implication of this algorithm is the identification of elements that construct 
quasipositive braids. Therefore, we believe that this algorithm may help to solve the 
quasipositivity problem in the braid group. 

\begin{\bib}{10}
\bibitem{Artin} Artin, E., {Theory of braids}, {\it Ann. Math.} {\bf 48} (1947), 101-126.
\bibitem{HT} Ben-Itzhak, T., Kaplan, S. and Teicher M., {Identifying Half-Twists Powers and 
Computing its Root}, preprint.
\bibitem{NEW} Birman, J.S., Ko, K.H. and Lee, S.J., {A new approach to the word and conjugacy 
problems in the braid groups}, {\it Adv. Math.} {\bf 139} (1998), 322-353.
\bibitem{Deh1} Dehornoy, P., {From large cardinals to braids via distributive algebra}, {\it J. 
Knot Theory \& Ramifications} {\bf 4(1)} (1995), 33-79.
\bibitem{Deh2} Dehornoy, P., {A fast method for comparing braids}, {\it Adv. Math.} {\bf 125(2)} 
(1997), 200--235. 
\bibitem{POSBR} Elrifai, E.A. and Morton, H.R., {Algorithms for positive braids}, {\it Quart. J. 
Math. Oxford Ser.} (2) {\bf 45} (1994), 479-497.
\bibitem{GAR} Garside, F.A., {The braid group and other groups}, {\it Quart. J. Math. Oxford Ser.} 
(2) {\bf 78} (1969), 235-254.
\bibitem{EFF} Jacquemard, A., {About the effective classification of conjugacy classes of braids}, 
{\it J. Pure. Appl. Alg.} {\bf 63} (1990), 161-169.
\bibitem{BAND} Kang, E.S., Ko, K.H. and Lee, S.J., {Band-generator presentation for the 4-braid 
group}, {\it Top. Appl.} {\bf 78} (1997), 39-60.
\bibitem{KUL} Kulikov, V.S. and Teicher, M., {Braid monodromy factorization and diffeomorphism 
types}, {\it Izvestiya: Mathematics} {\bf 64:2} (2000), 89-120. 
\bibitem{MILLER} Miller G.L., {Riemann's hypothesis and tests for primality}, {\it J. Computer and 
System Sciences} {\bf 13(3)} (1976), 300-317.
\bibitem{BGTI} Moishezon, B. and Teicher, M., {Braid group techniques in complex geometry I, Line 
arrangements in $\C \PP^2$}, {\it Contemporary Math.} {\bf 78} (1988), 425-555.
\bibitem{RABIN} Rabin M.O., {Probabilistic algorithm for testing primality}, {\it J. Number Theory} 
{\bf 12} (1980), 128-138.
\end{\bib}

\end{document}